\documentclass{amsart}

\usepackage{epsfig,psfrag}
\usepackage{latexsym}
\theoremstyle{plain}
\newtheorem{theorem}{Theorem}[section]

\newtheorem{prop}[theorem]{Proposition}

\newtheorem{alg}{Algorithm}

\theoremstyle{definition}
\newtheorem{defn}{Definition}[section]

\theoremstyle{remark}
\newtheorem{rem}{Remark}[section]
\newtheorem*{notation}{Notation}

\usepackage[all]{xy}

\theoremstyle{remark}

\usepackage{amssymb}
\usepackage{amsthm}
\usepackage{amsmath}
\usepackage{amsxtra}

\DeclareMathOperator{\chrt}{Char}

\DeclareMathOperator{\GL}{GL}
\DeclareMathOperator{\id}{Id}

\newcommand{\Q}{\mathbb Q}

\newcommand{\C}{\mathbb C}

\newcommand{\Z}{\mathbb Z}

{.tfm}
{.tfm}

\begin{document}
\title{The Dwork Family and Hypergeometric Functions}
\author{Adriana Salerno}

\maketitle

\begin{abstract}
 In his work studying the Zeta functions of families of hypersurfaces, Dwork came upon a one-parameter family of hypersurfaces (now known as \emph{the} Dwork family). These examples were not only useful to Dwork in his study of his deformation theory for computing Zeta functions of families, but they have also proven to be extremely useful to physicists working in mirror symmetry. A startling result is that these families are very closely linked to hypergeometric functions. This phenomenon was carefully studied by Dwork and Candelas, de la Ossa, and Rodr\'{i}guez-Villegas in a few special cases. Dwork, Candelas, et.al. observed that, for these families, the differential equation associated to the Gauss-Manin connection is in fact hypergeometric. We have developed a computer algorithm, implemented in Pari-GP, which can check this result for more cases by computing the Gauss-Manin connection and the parameters of the hypergeometric differential equation.

\end{abstract}

\section{Introduction}

In his work studying the Zeta functions of families of hypersurfaces, Dwork came upon a one-parameter family of hypersurfaces in $\mathbb{P}^{n-1}$ (now known as the \emph{Dwork family}), defined by:

$$X_{\lambda}:x_1^n+\cdots +x_n^n-n\lambda x_1\cdots x_n=0.$$

These examples were not only useful to Dwork in his study of his deformation theory for computing Zeta functions of families, but they have also proven to be extremely useful to physicists working in mirror symmetry (c.f. \cite{candelas} ). 

A startling result is that these families are very closely linked to hypergeometric functions. This phenomenon was carefully studied by Dwork in the cases where $n=3,4$ (see for example \cite{dwork:padic}) and for $n=5$ by Candelas, de la Ossa, and Rodr\'{i}guez-Villegas  (\cite{frv:dlo}, \cite{frv:candelas2}). 

In particular, Dwork's work gives a construction of modules isomorphic to the middle (relative) deRham cohomology, equipped with an integrable connection which is equivalent to the Gauss-Manin connection. By the work of Katz and Oda  \cite{katz:oda}, we know that this connection is essentially differentiation of cohomology classes with respect to the parameter. Differentiating each basis element in the module gives us a first-order system of differential equations. Dwork, Candelas et. al. noticed, in the examples they studied, that the differential equations obtained through this method are hypergeometric differential equations.  

In recent years, Dwork's ideas have been generalized to compute Zeta functions using $p$-adic and $\ell$-adic cohomology. In studying the Zeta function using $\ell$-adic cohomology, Katz proved that there was a link between more general monomial deformations of Fermat hypersurfaces (of which the Dwork family is an example) and hypergeometric sheaves \cite{katz:dwork}. Rojas-Leon and Wan, independently from Katz, implemented the same approach to compute Zeta-functions \cite{rojaswan}. In \cite{kloost}, Kloosterman shows that the $p$-adic Picard-Fuchs equation associated with the Dwork family is hypergeometric.

In this paper, we follow the direct approach originally used by Dwork and then Candelas et al, rather than the rigid cohomology approach.  That is, we use Dwork's original construction of a module over $\mathbb{C}$ and algorithmic methods based on combinatorics and linear algebra. . We have developed computer algorithms, which have been implemented in Pari-GP \cite{PARI2}.The GP scripts can be found in the Appendix.  One algorithm computes the matrix for the Gauss-Manin connection associated with the Dwork family, by blocks. Another uses a block of the connection matrix to compute the parameters of the associated hypergeometric differential equation (and in the process proving that this differential equation is hypergeometric).

\textbf{Acknowledgements.} As most of this work is based on the author's Ph.D. thesis \cite{salerno} , she would like to primarily thank her advisor Fernando Rodr\'{i}guez-Villegas for his guidance, support, and great ideas. Many people contributed to the progress of this thesis, and among them the author would like to thank Kiran Kedlaya, Frits Beukers, and Daqing Wan, and the Arizona Winter School for allowing her to meet these great mathematicians. Finally, the author would like to thank Michelle Manes, Bianca Viray, and ICERM for encouraging the final submission of this paper and Jonathan Webster for his computational number theory advice. 
\section{Background}

\subsection{Hypergeometric Functions}

\begin{defn}

Let, $A,B\in\Z$ and $\alpha_1,\dots,\alpha_A,\beta_1,\dots,\beta_B\in\Q$, with all of the $\beta_i\geq0$. The \textbf{generalized hypergeometric function} is defined as the series (taking $z\in\C$)

\begin{equation*}
{}_AF_B(\alpha_1,\dots,\alpha_A;\beta_1,\dots,\beta_B|z)=\sum_{k=0}^{\infty}\frac{(\alpha_1)_k\cdots(\alpha_A)_kz^k}{(\beta_1)_k\cdots(\beta_B)_kk!},
\end{equation*}

\end{defn} 

\noindent where we use the Pochhammer notation

\[(x)_k=x(x+1)\cdots(x+k-1)=\frac{\Gamma(x+k)}{\Gamma(x)}.\] 

The $\alpha_i$ will be referred to as ``numerator parameters'' and the $\beta_i$ as ``denominator parameters''.

Sometimes we will use the shortened notation $${}_AF_B(\alpha;\beta|z)={}_AF_B(\alpha_1,\dots,\alpha_A;\beta_1,\dots,\beta_B|z).$$

Let $\theta$ denote the operator $z\dfrac{d}{dz}$. The series ${}_AF_B(\alpha;\beta|z)$ satisfies the differential equation

\[\left\{\theta(\theta+\beta_1-1)\cdots(\theta+\beta_B-1)-z(\theta+\alpha_1)\cdots(\theta+\alpha_A)\right\}y=0.\]

Following the notation in \cite{beukers},

\[D(\alpha_1,\dots,\alpha_A;\beta_1,\dots,\beta_B)=\theta(\theta+\beta_1-1)\cdots(\theta+\beta_B-1)-z(\theta+\alpha_1)\cdots(\theta+\alpha_A).\]

If $A=B+1$ this is a Fuchsian differential equation with regular singularities at $z=0,1,\infty$ (in Section \ref{diffeq} we will review these definitions). We will focus only on hypergeometric functions with this property.

Notice that the parameters $\alpha_i,\beta_i$ completely characterize the hypergeometric function and its corresponding differential equation.

\subsection{Hypergeometric Groups}

In Section \ref{diffeq}, we will use a certain property of monodromy groups in order to relate the Gauss-Manin connection to hypergeometric functions. First, we need some definitions from \cite{beukers}.

Let $H$ be the fundamental group $\pi_1(\mathbb{P}^1\setminus\{0,1,\infty\},z_0)$ where $z_0$ is some fixed base point, for example $z_0=\frac{1}{2}$. Then clearly $H$ is generated by $g_0, g_1, g_{\infty}$ with the relation $g_{\infty}g_1g_0=1$, as pictured below.

\begin{figure}[h]
{\psfrag{z0}{$z_0$}
       \psfrag{0}{0}
       \psfrag{1}{1}
       \psfrag{g0}{$g_0$}
       \psfrag{g1}{$g_1$}
       \psfrag{ginf}{$g_{\infty}$}
\centerline{\includegraphics[width=20pc]{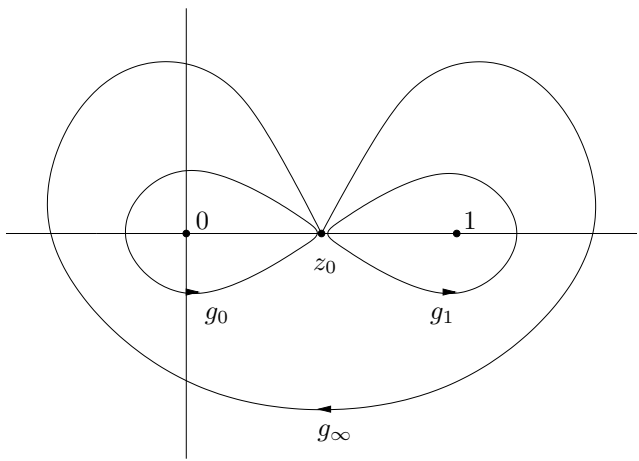}}}
\caption{The generators of $\pi_1(\mathbb{P}^1\setminus\{0,1,\infty\},z_0)$}
\end{figure}

Recall that the differential equation for a hypergeometric function of the form ${}_nF_{n-1}(\alpha;\beta|z)$ is Fuchsian with regular singular points $0,1,\infty$. Around a regular point, for example $z_0=\frac{1}{2}$, there are $n$ linearly independent analytic solutions with a non-zero radius of convergence. Let $A,B,C \in \GL(V)$ be determined by analytic continuation of solutions along the generators of  $\pi_1(\mathbb{P}^1\setminus\{0,1,\infty\},z_0)$, so that

\[\begin{array}{ccc}
    A&\leftrightarrow& g_{\infty}\\
    B&\leftrightarrow&g_0\\
    C&\leftrightarrow&g_1\\
    \end{array}\]

The group $\Gamma\subset \GL(V)$ generated by $A,B,C$ with the relation $ACB=\id$ is called the \emph{monodromy group}, and the map

\[H\rightarrow \GL(V)\]
\[g_{\infty},g_0,g_1\mapsto A,B,C,\]

\noindent is a representation of $H$.

\begin{defn}
Let $V$ be a finite dimensional complex vector space. A linear map $g\in \GL(V)$ is called a \emph{reflection} if $g-\id$ has rank one. The determinant of a reflection is called the \emph{special eigenvalue} of $g$.
\end{defn}

\begin{defn}
Suppose $a_1,\dots, a_n, b_1,\dots, b_n\in \C^*$ with $a_j\neq b_k$ for all $j=1,\dots,n$. A \emph{hypergeometric group} with numerator parameters $a_1,\dots,a_n$ and denominator parameters $b_1,\dots, b_n$ is a subgroup of $\GL(n,\C)$ generated by elements $h_0,h_1,h_{\infty}\in \GL(n,\C)$ such that $h_{\infty}h_1h_0=\text{Id}$,

\[\det(z-h_{\infty})=\prod_{i=1}^{n}(z-a_j)\]
\[\det(z-h_0^{-1})=\prod_{j=1}^{n}(z-b_j),\]

\noindent and $h_1$ is a reflection in the sense of Definition 4.1.
\end{defn}

Then we have the following useful result.

\begin{prop}\cite{beukers}
Suppose $a_1,\dots,a_n,b_1,\dots,b_n\in\C^*$ with $a_j\neq b_k$ for all $j,k=1,\dots,n$ and assume $b_n=1$. Let $\alpha_1,\dots,\alpha_n,\beta_1,\dots,\beta_{n-1}\in\C$ be such that $a_j=e^{2\pi i\alpha_j}$ for $j=1,\dots,n$ and $b_k=e^{2\pi i\beta_k}$ for $k=1,\dots,n-1$ . Then the monodromy group of the hypergeometric equation  \[D(\alpha_1,\dots,\alpha_n;\beta_1,\dots,\beta_{n-1})y=0\]

\noindent is a hypergeometric group with parameters $a_1,\dots,a_n,b_1,\dots,b_n$.
\end{prop}

The most important consequence of this result is that if we have a hypergeometric group $\Gamma$ in the sense of the Definition 2.3, we can find a hypergeometric differential equation whose monodromy group is $\Gamma$.

\section{Dwork's construction}

Suppose $X$ is a hypersurface defined over $\mathbb{P}^{n}$, and so it is $n-2$ complex dimensional. The \emph{middle} deRham cohomology will be the $(n-2)$-th cohomology. It is a classical result by Lefschetz that the $i$-th deRham cohomology of a projective hypersurface of dimension $n-2$ is identical to the $i$-th deRham cohomology of $\mathbb{P}^{n-1}$, for $i\neq n-2$. In other words, the middle cohomology is the only ``interesting'' one.

Recall the Dwork family of hypersurfaces defined by $$X_{\lambda}:x_1^n+\cdots +x_n^n-n\lambda x_1\cdots x_n=0.$$

It is not difficult to see that $X_{\lambda}$ is not smooth if and only if $\lambda$ is an $n$-th root of unity. Let $T=\C-\mu_n$. It follows that $X_{\lambda}$ is non singular for $\lambda\in T$. Dwork constructed modules over $\C$ isomorphic to the relative deRham cohomology $H^{n-2}_{dR}(X_{\lambda}/T)$, which are quite combinatorial in nature. 

To give Dwork's construction, we will use the notation established in \cite{katz2}. Let $F_{\lambda}(x_1,\dots,x_n):=x_1^n+\cdots +x_n^n-n\lambda x_1\cdots x_n$. Let $\mathfrak{L}$ be the free module (over $\C(\lambda)$) generated by the monomials $$x_1^{w_1}\cdots x_n^{w_n}=x^w,$$ with all the $w_i\geq0$ and $\sum_{i=1}^nw_i\equiv 0\bmod n$.

 Let $\mathfrak{L}^S$ be the submodule generated by monomials $x^w$ with all $w_i\geq1$. Let $D_i$ be the mapping defined by

 \[D_i:\mathfrak{L}\rightarrow\mathfrak{L},\hspace{1cm} D_i(x^w)=w_ix^w+x_i\frac{\partial F_{\lambda}}{\partial x_i}x^w.\]
 
  So $D_i(x^w)=w_ix^w+nx_i^nx^w-n\lambda x_1\cdots x_nx^w.$

 We now define the $\mathbb{C}(\lambda)$-module

 \[\mathcal{W}=\mathfrak{L}^S/\mathfrak{L}^S\bigcap\left(\sum_{i=1}^nD_i\mathfrak{L}\right).\]

 Thus, $\mathcal{W}$ is a vector bundle over $T$. 
 
 \begin{notation}
We will frequently represent a monomial $x^w$ by its exponent $w$.
\end{notation}

 \begin{alg}[Reduction Algorithm]\label{reduction}
 We define an algorithm on $n$-tuples $w$ using the relations given by the quotient described above. The input for this algorithm is $c* w=c*(w_1, w_2, \dots, w_n)$, where $c$ is an element of $\C(\lambda)$, and $w$ is an $n$-tuple representing a monomial in $\mathcal{W}$ where $w_i\geq n$ for some $i$. The output is a list that represents how to write $cw$ as a linear combination of monomials for which all exponents are less than $n$, thus ``reducing" $w$. 
 \begin{enumerate}
 \item Initialize two empty lists, $L$ and $M$. 
  \item Let $i$ be the first entry such that  $w_i\geq n$. Define $u=c\lambda* (w_1+1,w_2+1,\dots, w_i-n+1,\dots,w_n+1)$ and $v=c\frac{n-w_i}{n}* (w_1,\dots,w_i-n,\dots, w_n)$. If $\frac{n-w_i}{n}\neq0$:
 \begin{enumerate}
 \item If $v_i<n$ for all $i$, append $v$ to list $L$, unless $(w_1,\dots,w_i-n,\dots, w_n)$ is already in the list, then add $c\frac{n-w_i}{n}$ to the existing coefficient. Go to step 3.
 \item If $v_i\geq n$ for some $i$, append $v$ to list $M$, unless $(w_1,\dots,w_i-n,\dots, w_n)$ is already in the list, then add $c\frac{n-w_i}{n}$ to the existing coefficient. Go to step 3.
 \end{enumerate}
 \item If $u_i<n$ for all $i$, append $u$ to $L$ and proceed to step 4 (here if $u$ is already in $L$ we add coefficients as before). Otherwise, proceed to step 5.
 \item If the list $M$ is empty,  \textbf{output} the list $L$. Otherwise, proceed to step 8.
  \item If $u=w$, proceed to step 6. Otherwise, proceed to step 7. 
  \item If the list $M$ is empty, \textbf{output} the list $L$ but with all coefficients divided by $1-\lambda^n$. Otherwise, proceed to step 8.
  \item Append $u$ to $M$ and proceed to step 8.
  \item Let $a$ be the first element in $M$. Apply step 2 to $a$ instead of $w$ and remove $a$ from $M$. 
   \end{enumerate}
 \end{alg}  
 
 \begin{proof}
 Notice that we stop the algorithm whenever the list $M$ is empty.  The list $M$ contains monomials with entries that are greater than or equal to $n$.  To prove the algorithm terminates we need to prove that using this reduction process we can always empty the list.
 
 Given any starting $n$-tuple $w$, the worst that can happen is that every entry is greater than or equal to $n$. In this case, the first step would be to set $u=c\lambda* (w_1-n+1,w_2+1,\dots,w_n+1)$ and $v=c\frac{n-w_1}{n}* (w_1-n,\dots, w_n)$ . At most two monomials get added to $M$ ($v$ is added only if its coefficient is nonzero). We remove the first monomial (suppose this is $u$) and apply the reduction again. So $u_2=c\lambda^2*(w_1-n+2,w_2-n+2,\dots,w_n+2)$ and $v_2=c\lambda\frac{n-w_2}{n}* (w_1-n+1,w_2-n+1, \dots, w_n+1)$. Once more, at most two monomials get added to $M$. But now notice that the next item on the list would be $v$. Applying the reduction to this monomial would give $u_3=c\lambda\frac{n-w_1}{n}* (w_1-n+1,w_2-n+1,\dots, w_n+1)$ and $v_3=c\frac{n-w_1}{n}\frac{n-w_2}{n}* (w_1-n,w_2-n, \dots, w_n)$. Notice that $u_3=kv_2$, where $k$ is some constant, and thus for each pair of monomials we remove we are only really adding, at most, the number of monomials we had before plus one. 
 
 One can visualize the process in a diagram as follows:

{\Small{\xymatrix@C-2pc{ (w_1,w_2,\dots, w_n)\ar[d]_{(n-w_1)/n}\ar[dr]^{\lambda}& & \\
    (w_1-n,w_2,\dots,w_n)\ar[d]_{(n-w_2)/n}\ar[dr]^\lambda& (w_1-n+1,w_2+1,\dots,w_n+1)\ar[d]_{(n-w_2)/n}\ar[dr]^{\lambda} & \\
 (w_1-n,w_2-n,\dots,w_n)  & (w_1-n+1,w_2-n+1,\dots,w_n)&  (w_1-n+2,w_2-n+2,\dots,w_n+2) \\
    }}}

Notice the down arrows indicate subtracting $n$ from the first position that is greater than $n$, and this process clearly terminates. Also notice that the ``down and right" diagonals involve subtracting $n$ from the first entry larger than $n$ and then adding 1 to all entries. If this process does not terminate in a monomial whose entries are all less than $n$ then at the $n$-th reduction step we will find the monomial $(w_1-n+n,w_2-n+n,\dots, w_n-n+n)=w$, and as long as we can empty the list $M$, the algorithm terminates. But since all of the monomials are accounted for by following the ``down and right" diagonal and the down arrows, we are certain we can empty $M$. 
 \end{proof}

The following proposition is a direct consequence of  Algorithm \ref{reduction}. 

 \begin{prop}
 $\mathcal{W}$ is generated over $\C(\lambda)$ by the set of monomials

\[\mathcal{B}=\left\{x_1^{w_1}\cdots x_n^{w_n}=x^w| 1\leq w_i\leq n-1, \sum w_i\equiv0\bmod n\right\}.\]
In particular, $\mathcal{W}$ has dimension $$(n-1)^{(n-1)}-(n-1)^{(n-2)}+(n-1)^{(n-3)}-\cdots\pm(n-1).$$

 \end{prop}

The vector bundle $\mathcal{W}$ is equipped with an integrable connection $\nabla$ defined by

\[\nabla(f(\lambda)x^w)=\frac{\partial}{\partial\lambda}f(\lambda)x^w+f(\lambda)\frac{\partial}{\partial\lambda}F_{\lambda}x^w.\]

In particular, for monic monomials,

\[\nabla(x^w)=\frac{\partial}{\partial\lambda}F_{\lambda}x^w= -n x\cdot x^w.\]


Katz, in \cite{katz2}, proved the following useful lemmas. 

\begin{theorem}[The Comparison Theorem]
Let $w_0=\frac{1}{n}\sum_{i=1}^nw_i$. There is a $T$-linear map $\mathcal{R}:\mathfrak{L}^S\rightarrow H_{dR}^{n-1}(\mathbb{P}^n-X_{\lambda}/T)$ given by

\[\mathcal{R}:x^w\mapsto (-1)^{w_0}(w_0-1)!\frac{x^w}{F_{\lambda}^{w_0}}\frac{d(x_1/x_n)}{x_1/x_n}\wedge\cdots\wedge\frac{d(x_{n-1}/x_n)}{x_{n-1}/x_n}.\]
\end{theorem}

By this theorem and the residue map (cf. \cite{griffiths}) we have an isomorphism from $H_{dR}^{n-1}(\mathbb{P}^n-X_{\lambda}/T)$ to  $H_{dR}^{n-2}(X_{\lambda}/T)$.

And we have the following:
\begin{theorem}
The map $\Theta$ induces, by passage to quotients, an isomorphism

\[\Theta:\mathcal{W}\stackrel{\sim}{\rightarrow}H_{dR}^{n-2}(X_{\lambda}/T),\]

\noindent which is compatible with the connection.
\end{theorem}

Hence the space $\mathcal{W}$ obtained through Dwork's construction is isomorphic to the middle (relative) deRham cohomology.

It is also important to note that $\Theta$ transforms $\nabla$ into the Gauss-Manin connection.

\section{Computing the connection matrix}

Let $\mu_n^n$ denote the group of $n$-tuples of $n$-th roots of unity, and $\Delta$ denote the diagonal elements. 

The character group of $\mu_n^n/\Delta$ is in one-to-one correspondence with the set 

\[W=\{(w_1,\dots,w_n)|0\leq w_i<n,\sum w_i\equiv0\bmod n\},\]

\noindent where

\[\chi_w(\xi):=\chi(\xi^w),\hspace{1cm}\xi^w=\xi_1^{w_1}\cdots\xi_n^{w_n}\]

\noindent and $\chi$ is a fixed primitive character of $\mu_n$.

Let $$G=\{\xi\in\mu_n^n|\xi_1\cdots\xi_n=1\}/\Delta.$$ The characters $\chi_w$ of $\mu_n^n/\Delta$ which act trivially on $G$ are precisely powers of $\chi_{\overline{1}}$, where $\overline{1}=(1,1,\dots, 1)$. Thus, $\chrt(G)$, the character group of $G$, corresponds to equivalence classes of $w$ in $W$, where $w'\sim w$ if $w-w'$ is a multiple (mod n) of $\overline{1}$. 
 
The varieties $X_{\lambda}$ allow a faithful action of the group $G$ by $\xi=(\xi_1,\dots,\xi_n)$ taking the point $(x_1,\dots,x_n)$ to $(\xi_1x_1,\dots,\xi_nx_n)$. Using this action, we get that $\mathcal{W}$ splits into eigenspaces as follows.

\begin{prop}

The action of $G$ on a fiber $\mathcal{W}$, gives

 \[\mathcal{W}=\bigoplus_{\chi\in \chrt(G)} \mathcal{W}(\chi),\]

\noindent  where $\mathcal{W}(\chi)$ is an eigenspace with basis

\[\{w,w+\overline{1}(\bmod n),\dots,w+\overline{n-1}(\bmod n)\},\]

\noindent but we exclude adding any vector $\overline{m}$ such that $m+w_i\equiv0\bmod n$ for some $i$.
\end{prop}

 To understand $\nabla$'s effect on $\mathcal{W}$, it is sufficient to know what it does to elements in the basis $\mathcal{B}$. From the definition of $\nabla$ we see that

\[\nabla(x^w)=\frac{\partial}{\partial\lambda}F_{\lambda}x^w= -n x^{w+\overline{1}}\]

\noindent where $w+\overline{1}=(w_1+1,\dots,w_n+1)$.

   Thus, the connection commutes with the action of $G$, so the proposition implies that $\nabla$ preserves eigenspaces. We want to compute the connection matrix $\nabla$. Because of the way in which $\nabla$ preserves eigenspaces, the connection matrix will have blocks on its diagonal for each eigenspace.

The main idea of the following algorithm is to use the reduction algorithm described earlier on $\nabla(x^w)$ where $w$ runs through the basis of an eigenspace. 

\begin{alg}[Computing a block of the connection matrix]\label{connection}
This algorithm takes any vector of integers as an input and outputs a matrix that is the block of the connection matrix that corresponds to that vector's eigenspace generators.
\begin{enumerate}
\item Create a basis of the  eigenspace related to $w$ by computing $\mathcal{B}=\{w,w+\overline{1}(\bmod n),\dots,w+\overline{n-1}(\bmod n)\}=\{v_1,\dots,v_k\}$, where we omit monomials which have entries equal to 0 mod n (so $k$ may or may not equal $n$).
\item Create $M$, a $k\times k$ matrix of zeros. 
\item Let $i=1$.
\item If $i= k+1$, \textbf{output} $M$. Otherwise, take the monomial $v_i$ in $\mathcal{B}$ and compute its derivative, that is $\nabla(v_i)=-n*(v_i+\overline1)$.
\item If  all the entries in $\nabla(v_i)$ are less than $n$, then add $-n$ to the $(i,i+1)$ position in $M$, add 1 to $i$, and go back to step 5. If not, proceed to the next step.
\item Reduce $\nabla(v_i)$ according to Algorithm \ref{reduction}. From the algorithm, we have $\nabla(v_i)$ as a linear combination of elements of $\mathcal{B}$, so  that $\nabla(v_i)=\alpha_1v_1+\cdots\alpha_kv_k$. Replace the $i$-th row of $M$ by  $(\alpha_1,\dots,\alpha_k)$. Add 1 to $i$ and go to Step 5. 
\end{enumerate}

\end{alg}


\section{The differential equation associated to the connection}\label{diffeq}

In this section, we will show that the differential equation associated to the connection $\nabla$ is a hypergeometric differential equation. We have developed an algorithm which outputs the parameters $\alpha, \beta$ given each block representative. 

We will first establish some notation and definitions. For more details and proofs, see for example \cite{beukers:diff}, \cite{CL}, or \cite{ince}.  

Consider the $n$th order equation

\begin{equation}\label{ode}
\sum_{m=0}^na_{n-m}(z)y^{(m)}=0, \hspace{1cm} (a_0(z)\equiv1)
\end{equation}

\noindent where the $a_k(z)$ are single-valued and analytic in a punctured neighborhood of a point $z_0$. Recall that if any of the $a_k$ have a singularity at $z_0$, then $z_0$ is called a \emph{singular point} for (\ref{ode}), otherwise it is called an analytic point. We say $z_0$ is a \emph{regular singular point} if

\begin{equation}\label{ode2} a_k(z)=(z-z_0)^{-k}b_k(z), \hspace{1cm} (k=1,\dots,n),\end{equation}

\noindent where $b_k$ is analytic at $z_0$.

A \emph{system of n first order equations} over $\C(z)$ has the form
\begin{equation}\label{system}
y'= Ay
\end{equation}

\noindent in the unknown column vector $y=(y_1,\dots,y_n)^T$ and where $A$ is an $n\times n$-matrix with entries in $\C(z)$. The entries are assumed to be single-valued and analytic at a neighborhood of a point $z_0$, and will at most have a pole at that point.

If $A$ has a singularity at $z_0$, then $z_0$ is a \emph{singular point} for the system (\ref{system}). $z_0$ is a \emph{regular singular point} if

\begin{equation*}
A(z)=(z-z_0)^{-1}\tilde{A}(z)
\end{equation*}

\noindent where $\tilde{A}$ is analytic for a neighborhood of $z_0$ (including $z_0$), and $\tilde{A}(z_0)\neq0$.

To study the system at $\infty$, we change variables from $z$ to $1/\zeta$. The associated system is

\[\dfrac{d\tilde{y}}{d\zeta}=-\frac{\tilde{A}(\zeta)}{\zeta^2}\tilde{y}.\]

It is not difficult to see that a differential equation like (\ref{ode}) can be rewritten as a system by setting $y_1=y,y_2=y',\dots,y_n=y^{(n-1)}$. Notice that this means $y_1'=y_2, y_2'=y_3,\dots,y_{n-1}'=y_n$, and $y_n'$ is given by the differential equation. So the differential system is determined by a companion matrix, as follows:

\begin{equation}\label{system:matrix}
\frac{d}{dz}\left(\begin{array}{c}y_1\\y_2\\\vdots\\y_n\end{array}\right)=\left(\begin{array}{ccccc}
                                                                                0&1&0&\cdots&0\\
                                                                                0&0&1&\cdots&0\\
                                                                                \vdots&\vdots&\vdots&&\vdots\\
                                                                                -a_n&-a_{n-1}&-a_{n-2}&\cdots&-a_1
                                                                                \end{array}\right)
                                                                            \left(\begin{array}{c}y_1\\y_2\\\vdots\\y_n\end{array}\right).
\end{equation}

Now, we will explain how a system of first order differential equations arises from the connection. Recall  that, on a vector bundle, being equipped with a connection $\nabla$ is equivalent to being equipped with an action of $\dfrac{d}{d\lambda}$ (see for example \cite{kedlaya} for a detailed explanation). In short, we have a first-order system defined by

\[\dfrac{d}{d\lambda}y=A y,\]

\noindent where $A$ is actually the transpose of the matrix we found in the previous section. In fact, each block defines its own differential system. We will use the following fact:
\begin{theorem}[Cyclic Vector Lemma]
Any system of linear first order differential equations is equivalent to a system which comes from a differential equation.
\end{theorem}

A proof of the lemma  can be found in \cite{gfunctions}.

Basically, this theorem says that in the space of solutions of our differential system there is a cyclic vector, that is, a vector such that $v, Av, A^2v, \dots, A^{n-1}v$ is a basis. More specifically, this means that if $y=(y_1,y_2,\dots, y_n)^T$ is a solution for the system, we can find an equivalent system with solutions of the form $\hat{y}=\left(\hat{y}_1,\dfrac{d\hat{y}_1}{d\lambda},\dots,\dfrac{d^{n-1}\hat{y}_1}{d\lambda^{n-1}}\right)^T$. In fact, using the system, we can represent the derivatives $\hat{y}_1^{(k)}$ as a linear combination of $y_1,\dots,y_n$. This determines a change of basis matrix $S$ such that $Sy=\hat{y}$.

The vector $Sy=\hat{y}$ satisfies a differential system of the form

\[\dfrac{d}{d\lambda}\hat{y}=\left(SAS^{-1}+\dfrac{dS}{d\lambda}S^{-1}\right)\hat{y},\]

\noindent and this last system is the companion matrix to a higher order differential equation. In our situation, since the basis vectors are basically already derivatives of each other, any vector in the basis, for example $y_1$, is a cyclic vector, and so $S$ is easy to determine. Let $A_S=\left(SAS^{-1}+\dfrac{dS}{d\lambda}S^{-1}\right)$. 

Notice that, from the reduction algorithm, the entries in the connection matrix are polynomials in $\lambda$ or rational functions in $\lambda$. In the case that they are rational functions,  $1-\lambda^n$ is the only possible denominator. After doing the change of basis described above, we may get some powers of $\lambda$ in the denominator, but as the new system is a companion matrix, this will only happen in the last row of the matrix. 

At this point, we have a differential equation associated to the connection. But solving high order differential equations is not a simple task. Instead, we will show that these matrices are related to a hypergeometric group, which in turn gives us the defining parameters of the hypergeometric differential equation we want to find.

We have a way of changing from a differential equation to a system and viceversa. It is important to note that a regular singular point of (\ref{ode}) $z_0$ \emph{may not} be a regular singular point of the system associated with it. This happens only when the $a_k$ have at most simple poles at $z_0$. However, there is an equivalent first-order system with the property that if $z_0$ is a regular singular point of (\ref{ode}) then $z_0$ is a regular singular point of the system.

Suppose (\ref{ode}) has a regular singularity at $z_0$, and let $\phi$ be a solution of (\ref{ode}). Define $\hat\phi$ to be the vector with components $\phi_1,\dots,\phi_n$ by setting 

\[\phi_k=(z-z_0)^{k-1}\phi^{(k-1)},\hspace{1cm}(k=1,\dots,n).\]

Then, clearly,

\begin{eqnarray*}
(z-z_0)\phi_k'&=&(z-z_0)((z-z_0)^{k-1}\phi^{(k-1)})'\\
&=&(z-z_0)((k-1)(z-z_0)^{k-2}\phi^{(k-1)}+(z-z_0)^{k-1}\phi^{(k)})\\
&=&(k-1)(z-z_0)^{k-1}\phi^{(k-1)}+(z-z_0)^k\phi^{(k)}\\
&=&(k-1)\phi_k+\phi_{k+1}\hspace{1cm} (k=1\dots,n-1)
\end{eqnarray*}

And, finally,

\[(z-z_0)\phi_n'=(n-1)\phi_n-\sum_{m=1}^nb_{n-m+1}(z)\phi_m,\] where the $b_i$ are defined as in (\ref{ode2}). 

Therefore $\hat\phi$ is a solution of the linear system
\begin{equation}
y'=\hat{A}(z)y
\end{equation}

\noindent where $\hat A$ has the structure

\begin{equation*}
\hat{A}(z)=(z-z_0)^{-1}\left(\begin{array}{cccccc}
                                0&1&0&0&\cdots&0\\
                                0&1&1&0&\cdots&0\\
                                0&0&2&1&\cdots&0\\
                                0&0&0&3&\cdots&0\\
                                \vdots&\vdots&\vdots&\vdots&&\vdots\\
                                0&0&0&0&\cdots&1\\
                                -b_n&-b_{n-1}&-b_{n-2}&-b_{n-3}&\cdots&(n-1)-b_1
                                \end{array}\right)
\end{equation*}

This system clearly has a regular singularity at $z_0$.

In our situation, we want a hypergeometric differential system, so in particular we want a simple pole at $\lambda=0$. We can accomplish this by replacing the matrix $A_S$ by the matrix $\hat{A}$, which has a simple pole at 0. In particular, now we are certain that the only denominators in the last row are of the form $1-\lambda^n$. 

We will now change variables to $z=\lambda^n$. By the chain rule, we have:

\[\frac{d}{d\lambda}y=\frac{d}{dz}y\frac{d}{d\lambda}z=\frac{d}{dz}y*n\lambda^{(n-1)}\]

So our system representing the derivative with respect to $z$ is $y'=By$ where $B=\dfrac{1}{n\lambda^{(n-1)}}\hat{A}$. Since $\hat{A}$ has a simple pole at $\lambda=0$, replacing every instance of $\lambda^n$ by $z$ gives that $B$ has a simple pole at $z=0$ as well, but now we also have a simple pole at $z=1$.

There is an algorithm by Brieskorn which relates Gauss-Manin connections to monodromy group generators \cite{bries}. Let $A$ be the matrix representation of the connection. The algorithm uses the fact that if $A$ has a simple pole around a given point, i.e. can be written as

\[A=A_{-1}(z-z_0)^{-1}+A_0+A_1(z-z_0)+\cdots\]

\noindent then $R=e^{2\pi iA_{-1}}$ gives the monodromy around $z_0$.  

The process of changing the system matrix to $B$ ensures that we have a simple pole around zero, one, and infinity.  Since this last system has regular singular points at $0,1,\infty$ and no other singularities, it is Fuchsian, as we expected. 

Let $B_{-1}$ denote the residue around zero, and $\tilde{B}_{-1}$ denote the residue at $\infty$. Now define $h_0=e^{2\pi iB_{-1}}$ and $h_{\infty}=e^{2\pi i \tilde{B}_{-1}}$. Let $D$ be the residue of $B$ around $z=1$. The matrix $h_1=e^{2\pi iD}$ is clearly a reflection in the sense described by Beukers and Heckman. This is easily checked by noticing that there is only one row with denominators of the form $1-z$. Thus, the residue at one will necessarily have rank one, which implies that $h_1-\id$ has rank one. 
 
 We have shown that the matrices $h_{\infty}, h_1, h_0$ generate a hypergeometric group.

Recall that the monodromy group described by Beukers and Heckman in \cite{beukers} is generated by the monodromy matrices around zero, one, and infinity, and the parameters are determined by the eigenvalues of these matrices. In fact, we can get the parameters of the hypergeometric differential equation directly from $B_{-1}$, the residue around zero, and $\tilde{B}_{-1}$, where $\tilde{B}$ is the system at $\infty$. The eigenvalues of $B_{-1}$ will be the $\beta$'s and the eigenvalues of $\tilde{B}_{-1}$ will be the $\alpha$'s. Therefore, the monodromy group we have found corresponds to the differential equation
 
 \[D(\alpha;\beta|z)y=0,\] where the $\alpha$'s are the eigenvalues of $\tilde{B}_{-1}$ and the $\beta$'s are the eigenvalues of $B_{-1}$. 

Here is the algorithm we have just described.

\begin{alg}[Computing the parameters of the hypergeometric differential equation] This algorithm takes as input a monomial $w$ and computes the parameters of the hypergeometric differential equation associated to the eigenspace generated by $w$. The output is two lists, $\alpha$ and $\beta$, containing the parameters. 

\begin{enumerate}
\item Let $A$ be the $k\times k$ connection matrix block associated to $w$, computed with Algorithm \ref{connection}.
\item Compute the change of basis matrix, $S$, given by assuming $w$ as a cyclic vector.
\item Compute $A_S=\left(SAS^{-1}+\dfrac{dS}{d\lambda}S^{-1}\right)$.
\item Compute $\hat{A}$ by multiplying $a_{k,j}$ by $\lambda^{k-j+1}$ and adding $i-1$ in the $(i,i)$ position of $A_S$ for $i=1, \dots, k$ . This means that in the $(k,k)$ position we have $k-1+\lambda a_{k,k}$. 
\item Multiply $\hat{A}$ by $1/n$, call this $B$. 
\item Let $h_0=Res_{z=0}B$ and compute the eigenvalues. Let $\beta$ be the list of eigenvalues.
\item Let $h_{\infty}=Res_{z=\infty}B$ and compute the eigenvalues. Let $\alpha$ be the list of these eigenvalues.  
\item \textbf{Output:} $\alpha, \beta$. 

\end{enumerate}
\end{alg}

\begin{rem}
While computing some examples, we noticed that given a vector $(w_1,w_2,\dots,w_n)$, if we cancel out the numbers which it has in common with the list $(0,1,2,\dots,n-1)$, then $\alpha_i=\frac{w_j}{n}$ for each $w_j$ that survives the cancelation, and $\beta_i=\frac{k}{n}$ for each $k$ that survives in the second vector. This is exactly Katz's main result in the case of hypergeometric sheaves \cite{katz:dwork}. 
\end{rem}


\section{An Illustrative Example}

Suppose $n=6$. By Algorithm \ref{reduction} we have that

\[\mathcal{W}=\left\langle x^w|1\leq w_i\leq 5, \sum w_i\equiv0\bmod 6 \right\rangle,\]

\noindent and $\mathcal{W}$ has dimension $5^5-5^4+5^3-5^2+5=2605$.

 Take $w=(1,1,1,2,2,5)$. This belongs to the eigenspace generated by

\[\mathcal{B}_{(1,1,1,2,2,5)}=\{(1,1,1,2,2,5),(3,3,3,4,4,1),(4,4,4,5,5,2)\}.\]

Here is an example of the algorithm for computing the block in the matrix representation of $\nabla$ for $n=6$ and the eigenspace corresponding to the monomial $(1,1,1,2,2,5)\sim x_1x_2x_3x_4^2x_5^2x_6^5$, with basis denoted earlier by $\mathcal{B}_{(1,1,1,2,2,5)}$. I will denote this block by $\nabla_{\mathcal{B}_{(1,1,1,2,2,5)}}$.

\begin{enumerate}
\item Apply $\nabla(1,1,1,2,2,5)=-6(2,2,2,3,3,6)$. Using the relations we can write this last monomial in terms of the monomials in $\mathcal{B}_{(1,1,1,2,2,5)}$. We can represent the process of Algorithm \ref{reduction} graphically, as shown below:

    \xymatrix@C-2.5pc{&(2,2,2,3,3,6)\ar[dl]_{0}\ar[dr]^{\lambda}& \\
     (2,2,2,3,3,0)& &   (3,3,3,4,4,1)\\}

    \noindent This means that $(2,2,2,3,3,6)=\lambda(3,3,3,4,4,1)+ 0\cdot(2,2,2,3,3,0)$, which is a monomial in $\mathcal{B}_{(1,1,1,2,2,5)}$. 
    
    Thus, in the matrix representation of $\nabla_{\mathcal{B}_{(1,1,1,2,2,5)}}$, there will be a $-6\lambda$ as the $(2,1)$ entry.
\item We repeat this process for the next monomial in the basis, $(3,3,3,4,4,1)$. Applying the connection we get $\nabla(3,3,3,4,4,1)=-6(4,4,4,5,5,2)$. Since this monomial is already in $\mathcal{B}$ we write $-6$ in the $(3,2)$ position in the block matrix.
\item Take $\nabla(4,4,4,5,5,2)=-6(5,5,5,6,6,3)$. We have to do the reduction process again, represented below.

{\tiny{
\xymatrix@C-2.5pc{ &(5,5,5,6,6,3)\ar[dl]_{0}\ar[dr]^{\lambda}& & & & & & &\\
    (5,5,5,0,6,3)& & (6,6,6,1,7,4)\ar[dl]_{-1/6}\ar[dr]^{\lambda} & & & & & &\\
    &(6,6,6,1,1,4)\ar[dl]_{0}\ar[dr]^{\lambda}& & (7,7,7,2,2,5)\ar[dl]_{-1/6}\ar[dr]^{\lambda} & & & & &\\
    (0,6,6,1,1,4)& &(1,7,7,2,2,5)\ar[dl]_{-1/6}\ar[dr]^{\lambda}& &(2,8,8,3,3,6)\ar[dl]_{-2/6}\ar[dr]^{\lambda}& & & &\\
    & (1,1,7,2,2,5)\ar[dl]_{-1/6}\ar[dr]^{\lambda}& & (2,2,8,3,3,6)\ar[dl]_{-2/6}\ar[dr]^{\lambda}& &(3,3,9,4,4,7)\ar[dl]_{-3/6}\ar[dr]^{\lambda}& & \\
    (1,1,1,2,2,5)& &(2,2,2,3,3,6)\ar[dl]_{0}\ar[dr]^{\lambda}& &(3,3,3,4,4,7)\ar[dl]_{-1/6}\ar[dr]^{\lambda}& &(4,4,4,5,5,8)\ar[dl]_{-2/6}\ar[dr]^{\lambda}& \\
    & (2,2,2,3,3,0)& &(3,3,3,4,4,1)& &(4,4,4,5,5,2)& &(5,5,5,6,6,3) \\  }
}}

This is a bit harder to unravel than the other cases, but it works in exactly the same way. The diagram shows us that

\begin{eqnarray*}
(5,5,5,6,6,3)&=&-\frac{\lambda^2}{108}(1,1,1,2,2,5)+\frac{17\lambda^4}{36}(3,3,3,4,4,1)\\
&&-\frac{3\lambda^5}{2}(4,4,4,5,5,2)+\lambda^6(5,5,5,6,6,3).
\end{eqnarray*}

And solving for $(5,5,5,6,6,3)$ we get that

\begin{eqnarray*}
&&\nabla(4,4,4,5,5,2)=-6(5,5,5,6,6,3)\\
&&=-\frac{\lambda^2}{18(\lambda^6-1)}(1,1,1,2,2,5)+\frac{17\lambda^4}{6(\lambda^6-1)}(3,3,3,4,4,1)\\
&&-\frac{9\lambda^5}{\lambda^6-1}(4,4,4,5,5,2).
\end{eqnarray*}

\item Combining all of these steps, we can write $\nabla_{\mathcal{B}_{(1,1,1,2,2,5)}}$ as

\[\nabla_{\mathcal{B}_{(1,1,1,2,2,5)}}=\left(\begin{array}{ccc}
                0&0&-\dfrac{\lambda^2}{18(\lambda^6-1)}\\
                -6\lambda&0&\dfrac{17\lambda^4}{6(\lambda^6-1)}\\
                0&-6&-\dfrac{9\lambda^5}{\lambda^6-1}\\
                \end{array}
        \right).\]

\end{enumerate}

We will now describe the algorithm for finding the parameters of the differential equation for the same example monomial. The connection above gives us that differentiation with respect to $\lambda$ is equivalent to the differential system

\[\dfrac{d}{d\lambda} y=\left(\begin{array}{ccc}
                0&-6\lambda&0\\
                0&0&-6\\
                -\dfrac{\lambda^2}{18(\lambda^6-1)}&\dfrac{17\lambda^4}{6(\lambda^6-1)}&-\dfrac{9\lambda^5}{\lambda^6-1}\\
                \end{array}\right) y.\]

We compute the change of basis for the cyclic vector lemma, which is

\[S=\left(\begin{array}{ccc}
            1&0&0\\
            0&-6\lambda&0\\
            0&-6&36\lambda\\
            \end{array}\right).\]

And so we have the new system

\[\dfrac{d}{d\lambda}y=\left(SAS^{-1}+\dfrac{dS}{d\lambda}S^{-1}\right)y=\left(\begin{array}{ccc}
            0&1&0\\
            0&0&1\\
            \dfrac{2\lambda^3}{1-\lambda^6}&\dfrac{10\lambda^6-2}{\lambda^2(1-\lambda^6)}&\dfrac{7\lambda^6+2}{\lambda(1-\lambda^6)}\\
            \end{array}\right)y,\]

\noindent which, as we expected, is given by the companion matrix for a third order  differential equation.

This system is equivalent to

\[\dfrac{d}{d\lambda}y=\frac{1}{\lambda}\left(\begin{array}{ccc}
            0&1&0\\
            0&1&1\\
            \dfrac{2\lambda^6}{1-\lambda^6}&\dfrac{10\lambda^6-2}{1-\lambda^6}&2-\dfrac{7\lambda^6+2}{1-\lambda^6}\\
            \end{array}\right)y,\] which clearly has a simple pole at $\lambda=0$. 

Now we can change variables by setting $z=\lambda^6$. The change of variables leaves us with a system

\[\dfrac{d}{d\lambda}y=\frac{1}{z}\left(\begin{array}{ccc}
            0&1/6&0\\
            0&1/6&1/6\\
            \dfrac{z}{3(1-z)}&\dfrac{5z-1}{3(1-z)}&\dfrac{5z+4}{6(1-z)}\\
            \end{array}\right)y.\]

 The residue at zero is

\[A_{-1}=\left(\begin{array}{ccc}
            0&1/6&0\\
            0&1/6&1/6\\
            0&-{1}/{3}&2/3\\
            \end{array}\right),\]

\noindent which has eigenvalues $0,1/2,1/3$. 

Also, around infinity we have

\[\dfrac{d\tilde{y}}{d\zeta}=\frac{1}{\zeta}\left(\begin{array}{ccc}
            0&-1/6&0\\
            0&-1/6&-1/6\\
            \dfrac{1}{3(1-\zeta)}&\dfrac{5-\zeta}{3(1-\zeta)}&\dfrac{5+4\zeta}{6(1-\zeta)}\\
            \end{array}\right)\tilde{y},\]

\noindent which has residue (at $\zeta=0$) of

\[\tilde{A}_{-1}=\left(\begin{array}{ccc}
            0&-1/6&0\\
            0&-1/6&-1/6\\
            {1}/{3}&{5}/{3}&{5}/{6}\\
            \end{array}\right),\]

\noindent and thus yields the eigenvalues $1/3,1/6,1/6$.

            We have now found the parameters of the hypergeometric differential equation associated to this connection matrix block:
\[D\left(\frac{1}{6},\frac{1}{6},\frac{1}{3};\frac{1}{2},\frac{2}{3}\right)y=0.\]

To sum it up, the block of the matrix $\nabla$ corresponding to the eigenspace of $(1,1,1,2,2,5)$ gives rise to the hypergeometric differential equation which has

\[{}_3F_2\left(\left.\frac{1}{6},\frac{1}{6},\frac{1}{3};\frac{1}{2},\frac{2}{3}\right|z\right)\]

\noindent as its holomorphic solution around 0.

Table 1 below shows some numerical examples for $n=6$. 
\begin{table}[ht]\label{table}
\renewcommand{\arraystretch}{2}
\caption{Parameters for $n=6$}
\centering
\begin{tabular}{|c|c|c|}
  \hline

  Monomial & $\alpha_i$ & $\beta_i$ \\

  \hline

  $[1, 1, 1, 1, 1, 1]$ & $\left[\dfrac{1}{6}, \dfrac{1}{6}, \dfrac{1}{6}, \dfrac{1}{6}, \dfrac{1}{6}\right]$ & $\left[ \dfrac{1}{2}, \dfrac{2}{3}, \dfrac{5}{6}, \dfrac{1}{3}\right]$ \\

  \hline

  $[5, 3, 1, 1, 1, 1]$ &$\left[\dfrac{1}{6}, \dfrac{1}{6}, \dfrac{1}{6}\right]$& $\left[\dfrac{2}{3}, \dfrac{1}{3}\right]$ \\

  \hline

  $[4,4,1,1,1,1]$ & $\left[\dfrac{1}{6},\dfrac{1}{6},\dfrac{1}{6}, \dfrac{2}{3}\right]$ & $\left[\dfrac{1}{3},\dfrac{1}{2},\dfrac{5}{6}\right]$\\

  \hline

  $[5, 2, 2, 1, 1, 1]$  & $\left[\dfrac{1}{3}, \dfrac{1}{6}, \dfrac{1}{6}\right]$ & $\left[\dfrac{1}{2}, \dfrac{2}{3}\right]$ \\
  \hline
  $[4, 3, 2, 1, 1, 1]$ & $\left[\dfrac{1}{6}, \dfrac{1}{6}\right]$ & $\left[\dfrac{5}{6}\right]$ \\
  \hline
  $[3, 3, 3, 1, 1, 1]$     & $\left[\dfrac{1}{2}, \dfrac{1}{2}, \dfrac{1}{6}, \dfrac{1}{6}\right]$ & $\left[\dfrac{2}{3}, \dfrac{5}{6}, \dfrac{1}{3}\right]$ \\
  \hline
  $[4, 2, 2, 2, 1, 1]$  & $\left[\dfrac{1}{3}, \dfrac{1}{3}, \dfrac{1}{6}\right]$& $\left[\dfrac{1}{2}, \dfrac{5}{6}\right]$ \\
  \hline
  $[3, 3, 2, 2, 1, 1]$   & $\left[\dfrac{1}{2}, \dfrac{1}{3}, \dfrac{1}{6}\right]$ & $\left[\dfrac{2}{3}, \dfrac{5}{6}\right]$ \\ \hline

  $[3, 2, 2, 2, 2, 1]$  & $\left[\dfrac{1}{3}, \dfrac{1}{3}, \dfrac{1}{3}\right]$ & $\left[\dfrac{2}{3}, \dfrac{5}{6}\right]$ \\ \hline
  $[5,5,3,3,1,1]$ & $\left[\dfrac{1}{6},\dfrac{1}{2},\dfrac{5}{6}\right]$ & $\left[\dfrac{1}{3},\dfrac{2}{3}\right]$\\ \hline
  $[5, 5, 4, 2, 1, 1]$ & $\left[\dfrac{1}{6}, \dfrac{5}{6}\right]$ & $\left[\dfrac{1}{2}\right]$ \\ \hline
  $[5, 4, 4, 3, 1, 1]$ & $\left[\dfrac{1}{6}, \dfrac{2}{3}\right]$ & $\left[\dfrac{1}{3}\right]$ \\ \hline
  $[5, 4, 3, 3, 2, 1]$ & $\left[\dfrac{1}{2}\right]$ & $[]$ \\ \hline
  $[4, 4, 4, 3, 2, 1]$ & $\left[\dfrac{2}{3}, \dfrac{2}{3}\right]$ & $\left[\dfrac{5}{6}\right]$\\
  \hline
  \end{tabular}
\end{table}

\bigskip

\section*{Appendix\\ GP Scripts}

\subsection*{Computing the connection matrix}

This function counts the number of coordinates of a vector with entries bigger than or equal to $d$.
\begin{verbatim}
count(a)=
{	
	   local(t);
	   t=0;
	   for(k=1,length(a), if(a[k]>=length(a),t=t+1));
	   t
}
\end{verbatim}

This function (from \cite{frv:ent}) can tell if a given element is in a vector, and gives the ``position'' of the element.

\begin{verbatim}
memb(g,v)=for(k=1,length(v),if(g==v[k],return(k)));0
\end{verbatim}

The following function takes a vector $a$ and an integer $m$ (its coefficient) and subtracts one to all the entries and adds d to one
of them until it gets to 0, a vector in the basis, or the original vector. It saves the leftovers in a vector $v$. It is one of the two possible reductions coming from the relations on $\mathcal{W}$.

\begin{verbatim}
red1(a,m)=
{	
    local(j,l,b,h,t,d,s,v,u);
    h=m;
    b=0;
    t=vector(length(a));
    d=length(a);
    l=1;
    j=1;
    s=vector(length(a));
    v=vector(0);
    t=a;
    if(count(t)>0,until(count(t)==0 || t==a,
        for(k=1,length(t), if(t[k]<t[j], j=k));
        for(k=1,length(t),s[k]=t[k]-1);
        for(k=1,length(t), if(k==j,t[k]=s[k]+d, t[k]=s[k]));
        l=h*s[j]/(d*n);
        v=concat(v,[[l,s]]);
        h=h*1/n));
    [h,t]		
}


red1leftovers(a,m)=
{   local(j,l,b,h,t,d,s,v,u);
    h=m;
    b=0;
    t=vector(length(a));
    d=length(a);
    l=1;
    j=1;
    s=vector(length(a));
    v=vector(0);
    t=a;
    if(count(t)>0,until(count(t)==0 || t==a,
        for(k=1,length(t), if(t[k]<t[j], j=k));
        for(k=1,length(t),s[k]=t[k]-1);
        for(k=1,length(t), if(k==j,t[k]=s[k]+d, t[k]=s[k]));
        l=h*s[j]/(d*n);
        v=concat(v,[[l,s]]);
        h=h*1/n));
    v		
}
\end{verbatim}

Here is the other possible reduction. This one subtracts $n$ from one spot and adds one to everything afterwards. Saves leftovers in vector $v$.

\begin{verbatim}
red2(a,m)=
{	
    local(j,l,b,h,t,d,s,v,u);
    h=m;
    b=0;
    d=length(a);
    l=1;
    t=vector(length(a));
    v=vector(0);
    j=1;
    s=vector(length(a));
    t=a;
    if(count(t)>0,until(count(t)==0 || t==a,
            for(k=1,length(t), if(t[j]<t[k], j=k));
            for(k=1,length(t),if(k==j,s[k]=t[k]-d,s[k]=t[k]));
            for(k=1,length(t), t[k]=s[k]+1);
            l=h*(-s[j])/d;
            v=concat(v,[[l,s]]);
            h=h*n));
    [h,t]
}

red2leftovers(a,m)=
{	
    local(j,l,b,h,t,d,s,v,u);
    h=m;
    b=0;
	   d=length(a);
	   l=1;
	   t=vector(length(a));
	   v=vector(0);
	   j=1;
	   s=vector(length(a));
	   t=a;
	   if(count(t)>0,until(count(t)==0 || t==a,
		      	      for(k=1,length(t), if(t[j]<t[k], j=k));
			             for(k=1,length(t),if(k==j,s[k]=t[k]-d,s[k]=t[k]));
			             for(k=1,length(t), t[k]=s[k]+1);
			             l=h*(-s[j])/d;
			             v=concat(v,[[l,s]]);
			             h=h*n));
	       v
}
\end{verbatim}

Now we combine these two reductions and loop until we get monomials in  the basis of $\mathcal{W}$. The input of this function is a vector of any length and the output will be the ``linear combination'' of that vector in terms of the basis vectors (vectors with entries between 1 and the length).

\begin{verbatim}
reduction(a)=
{
local(d,b,c,u,v,w,uu, bb,j, t, s, g,r);
u=vector(0);
v=vector(0);
d=length(a);
j=1;
if(count(a)==d,  u=[red1(a,1)]; v=red1leftovers(a,1),
                    u=[red2(a,1)];v=red2leftovers(a,1));
for(k=1, 10^d,
    if(k<=length(v),
        if(count(v[k][2])==0,
            if(v[k][1]==0, ,b=0;
                for(i=1,length(u),
                    if(v[k][2]==u[i][2],
                        u[i][1]=u[i][1]+v[k][1],
                        b=b+1));
                if(b==length(u),u=concat(u,[v[k]]))),
            if(v[k][1]==0, ,
                uu=red2(v[k][2],v[k][1]);
                v=concat(v,red2leftovers(v[k][2],v[k][1]));
                b=0;
                for(i=1,length(u),
                    if(uu[2]==u[i][2],
                        u[i][1]=u[i][1]+uu[1],
                        b=b+1));
                    if(b==length(u),u=concat(u,[uu])))),
            break));
b=0;
for(k=1, length(u),
    if(u[k][2]==a,
        w=vector(length(u)-1);
        for(j=1,k-1,w[j]=[u[j][1]/(1-u[k][1]),u[j][2]]);
        for(j=k,length(u)-1,
            w[j]=[u[j+1][1]/(1-u[k][1]),u[j+1][2]]),
        b=b+1));
        if(b==length(u), r=u, r=w);		
    r
}
\end{verbatim}

The next step is to write the connection matrix from this, that is, write a function that gives the derivatives of each vector in terms of the basis. In fact, there is an easy way to write the derivative of any vector using the reduction function.

\begin{verbatim}
derivative(a)=
{  	
    local(d,t,w);
    d=length(a);
    t=vector(d);
    for(k=1,d,t[k]=a[k]+1);
    if(count(t)==0,[[-d,t]],
        w=reduction(t);
        for(k=1,length(w),
        w[k][1]=w[k][1]*(-d)); w)
}
\end{verbatim}

Given a basis vector, we can find all the other basis vectors that will be a basis for the same eigenspace.

\begin{verbatim}
orbit(a)=
{   local(l,m, c,ss);
    d=length(a);
    l=0;
    c=0;
    ss=1;
    for(k=1,d, for(t=1,d, if(k==d-a[t], l=l+1;break)));
    m=d-l;
    b=vector(m); for(k=1,m, b[k]=vector(d));
    b[1]=a;
    for(s=1,d-1,for(t=1,d, if(s==d-a[t], , c=c+1));
        if(c==d,ss=ss+1;for(t=1,d, b[ss][t]=(a[t]+s)%d));c=0);
    b;
}
\end{verbatim}

The following gives the matrix representation of the block of the Gauss-Manin connection associated to a particular basis vector (i.e., it gives a block of the whole matrix, which is related to the eigenspace related to this basis vector).

\begin{verbatim}
connectionmatrix(a)=
{   local(v,w,M);
    v=orbit(a);
    M=matrix(length(v),length(v));
    for(j=1,length(v),
        w=derivative(v[j]);
        for(k=1,length(w),M[memb(w[k][2],v),j]=w[k][1]));
    M=mattranspose(M);
    M
}
\end{verbatim}	

\subsection*{The algorithm to find the differential equation}

The following finds the derivative with respect to $\lambda$ of a vector with a coefficient. It is essentially the product rule.

\begin{verbatim}
derivn(a)=
{   local(b,z, ww, vv);
    b=deriv(a[1]);
    z=derivative(a[2]);
    ww=[[b,a[2]]];
    vv=vector(length(z));
    for(k=1,length(z), vv[k]=[a[1]*z[k][1],z[k][2]]);
    for(k=1,length(vv),
        if(vv[k][2]==a[2], ww[1][1]=ww[1][1]+vv[k][1],
            ww=concat(ww,[vv[k]])));
    ww
}
\end{verbatim}

We would like to have the derivative of a vector which is a linear combination of these monomials. This should use ideas like the function above. The first function finds the derivative of a vector (with a coefficient) in a prescribed basis determined by the orbit of $b$. The second does the same, but only outputs the vector of coordinates, without writing the basis down.

\begin{verbatim}
derivv(a,b)=
{   local(v,w,z);
    w=orbit(b);
    v=vector(length(w));
    for(i=1,length(w), v[i]=[0,w[i]]);
    for(k=1,length(a),
        z=derivn(a[k]);
        for(j=1,length(z),
            v[memb(z[j][2],w)][1]=v[memb(z[j][2],w)][1]+z[j][1]
        );
    );
    v
}

derivv2(a,b)=
{   local(v,w,z);
    w=orbit(b);
    v=vector(length(w));
    for(k=1,length(a),
        z=derivn(a[k]);
        for(j=1,length(z),
            v[memb(z[j][2],w)]=v[memb(z[j][2],w)]+z[j][1]
        );
    );
    v
}
\end{verbatim}

We want to change basis, and we need a matrix that changes from our basis obtained by using ``connection'' to a basis obtained from derivatives (i.e. we use the cyclic vector theorem to write the matrix as the companion matrix to a differential equation).

\begin{verbatim}
cob(a)=
{
    local(z, vv, uu, w);
    uu=orbit(a);
    r=vector(length(uu),k,if(k==1,1));
    vv=vector(length(uu), k, [r[k],uu[k]]);
    z=[r];
        for(k=1,length(uu)-1,
                w=derivv2(vv,a);
                vv=derivv(vv,a);
                z=concat(z,[w])
        );
    S=Mat(z~);
    S
}
\end{verbatim}

This function takes two input vectors, one is a basis and the other a vector indicating a linear combination of elements in this basis.

\begin{verbatim}
cobv(uu,r)=
{
    local(z, vv, w);
    vv=vector(length(uu), k, [r[k],uu[k]]);
    z=[r];
        for(k=1,length(uu)-1,
                w=derivv2(vv,a);
                vv=derivv(vv,a);
                z=concat(z,[w])
        );
    S=Mat(z~);
    S
}
\end{verbatim}

The following computes what a change of basis does to the system of differential equations, where we change from a basis found by using the connection function to a basis of all the derivatives of a specific vector.

\begin{verbatim}
cobsystem(A, S)=
{
	local(dS, C);
	dS=matrix(length(A), length(A), X, Y, deriv(S[X,Y]));
	C=S*A*1/S+dS*1/S;
    C
}
\end{verbatim}

Now, as seen in Section \ref{diffeq}, we can change this system into an equivalent one with a simple pole at 0. 
\begin{verbatim}
 regform(A)=
{
local(m,Areg);
m=length(A);
Areg=A;
for(k=1,m, Areg[m,k]=Areg[m,k]*n^(m-k+1));
for(k=1,m,
	   for(i=1,m, if (i==k, Areg[k,i]=Areg[k,i]+k-1)));
Areg
}
\end{verbatim}

We need to change variables to have it in terms of $z$ instead of $\lambda^d$.

\begin{verbatim}
varchange(A,d)=
{
for(k=1,length(A),
        for(i=1,length(A),
            A[k,i]=substpol(A[k,i],n^d,x)/d));
    A
}
\end{verbatim}

Finally we can compute the residue of the matrix at $z=0$.

\begin{verbatim}
residuezero(A)=
{
	for(k=1,length(A),
		for(i=1,length(A),
			A[k,i]=subst(A[k,i],x,0)));
	A
}		
\end{verbatim}

 We want a function that finds the rational roots of a polynomial with rational coefficients (because all of our characteristic polynomials are of that form and only have rational roots). We are using the rational roots theorem.

\begin{verbatim}
ratlroots(f)=
{
    local(p,q,z,vv,a,b,c,n,r,j,i);
    c=poldegree(f);
    vv=vector(c+1);
    r=vector(0);
    z=vector(0);
        for(k=1,c+1, vv[k]=polcoeff(f,k-1));
    n=denominator(vv);
    f=f*n;
    a=substpol(f,x,0);
        if(a==0, r=concat(r,0); f=f/x; a=substpol(f,x,0));
    b=pollead(f);
    p=concat(divisors(a),-divisors(a));
    q=concat(divisors(b),-divisors(b));
    for(k=1,length(p),
        for(i=1,length(q),
            if(memb(p[k]/q[i],z)==0, z=concat(z,p[k]/q[i]))));
    for(i=1,length(z),
        if(substpol(f,x,z[i])==0, r=concat(r,z[i]);
                        f=f/(x-z[i]);i=1));
    r
}
\end{verbatim}

We now want to combine all these steps to find the hypergeometric parameters given a vector (or monomial) in $\mathcal{W}$.

\begin{verbatim}
hypergcoeff(a)=
{
    local(A, B, S, Azero, Ainf,f,g, vv, uu, d, m,r,t);
    d=length(a);
    A=connectionmatrix(a);
    S=cob(a);
    A=cobsystem(A,S);
    A=regform(A);
    A=varchange(A,d);
    Azero=residuezero(A);
        f=charpoly(Azero);
    B=matrix(length(A),length(A));
    Ainf=matrix(length(A),length(A));
    for(k=1,length(A),
        for(i=1,length(A),
            B[k,i]=-substpol(A[k,i],x,1/y)));
    for(k=1,length(A),
        for(i=1,length(A),
            Ainf[k,i]=subst(B[k,i],y,0)));
        g=charpoly(Ainf);
    r=ratlroots(f);
    m=vector(length(r));
    for(k=1,length(r), m[k]=1-r[k]);		
\\	print("alphas  ",ratlroots(g)," betas ",m)
    t=[ratlroots(g),m];
    t
}
\end{verbatim}

\subsection*{Creating a hypergeometric table} To make a table like  Table \ref{table} we should be able to check the hypergeometric coefficients somewhat systematically. So we have to find a good way to generate basis vectors (or representatives up to permutations of the variables). We are certain that there are more efficient ways to do this, based on conversations with computational number theorists, but the focus of our project is to obtain the relationship with hypergeometric functions, which this code does. 
 
This function (also from \cite{frv:ent}) finds all the partitions of a number $m$.

\begin{verbatim}
part(m)=
{
  local(k,j,sm,sj,s, S = []);
  k = j = 1;
  sm = sj = vector(m+1);

  while(k,
      s = sm[k]+j;
      if (s > m,
      until(j <= m, j = sj[k]+1; k--);
      next);

      k++; sm[k]=s; sj[k]=j;

      if (s == m,
      S = concat(S, [vector(k-1,l, sj[k-l+1])])));

  S
}
\end{verbatim}

This function uses the previous one to find the partitions of a number $m$ of length $c$ and into numbers that are less than $c$.

\begin{verbatim}
part2(m,c)=
{
    local(v,w,t);
    v=part(m);
    w=vector(0);
    for(k=1,length(v),
        if(length(v[k])==c,
            t=count(v[k]);
            if(t==0,
                w=concat(w,[v[k]]))));
    w
}
\end{verbatim}

Now we can put this together to get representatives of the basis. We don't have strict representatives, but at least we eliminate the cases in which none of the entries are equal to one, because those obviously are in an eigenspace with a vector with entries equal to one.

\begin{verbatim}
basisreps(m)=
{
    local(v,w);
    v=vector(0);
    for(k=1,ceil((m-1)/2),
        w=part2(m*k,m);
        for(j=1, length(w),
            if(memb(1,w[j])==0, ,v=concat(v,[w[j]]))));
    v
}
\end{verbatim}

We now want to be able to output a table with the basis vectors and the hypergeometric parameters associated to them given a number $d$. It should turn out to be the list of numbers in the vector that remain after canceling out with the list of numbers between 0 and $d$, in accordance with Katz's results \cite{katz:dwork}.

\begin{verbatim}
hypergtable(d)=
{
    local(v,u);
    v=basisreps(d);
    for(k=1,length(v), u=hypergcoeff(v[k]);
    print(v[k],"       ", "alphas       ", u[1],
                         "     betas    ", u[2]))
}	
\end{verbatim}


\begin{thebibliography}{10}

\bibitem{kedlaya:roe}
	T.A.~Abbott; K.~Kedlaya; and D.~Roe.
	\newblock Bounding {P}icard numbers of surfaces using $p$-adic cohomology.
	\newblock{\em Arithmetic, Geometry and Coding Theory (AGCT 2005)}, arXiv:math/0601508v2, to appear.
	


\bibitem{beukers}
F.~Beukers and G.~Heckman.
\newblock Monodromy for the hypergeometric function ${}_n{F}_{n-1}$.
\newblock {\em Invent. Math.}, 95:325--354, 1989.

\bibitem{beukers:diff}
Frits Beukers.
\newblock Ordinary linear differential equations.
\newblock {\em Course Lecture Notes}.

\bibitem{bries}
Egbert Brieskorn.
\newblock Die monodromie der isolierten singularitaten von hyperflachen.
\newblock {\em Manuscripta Math.}, 2:103--161, 1970.

\bibitem{candelas}
Philip Candelas and Xenia de~la Ossa.
\newblock The zeta-function of a $p$-adic manifold, {D}work theory for
  physicists.
\newblock {\em arxiv:0705.2056v1}, 2008.



\bibitem{CL}
E.A. Coddington and N.~Levinson.
\newblock {\em Theory of ordinary differential equations}.
\newblock McGraw-Hill, 1955.

\bibitem{frv:dlo}
P.~Candelas;~X. de~la Ossa; and F.~Rodr\'{i}guez-Villegas.
\newblock {C}alabi-{Y}au manifolds over finite fields {I}.
\newblock {\em http://xxx.lanl.gov/abs/hep-th/0012233}.

\bibitem{frv:candelas2}
P.~Candelas;~X. de~la Ossa; F. Rodr\'{i}guez-Villegas.
\newblock {C}alabi-{Y}au manifolds over finite fields {II}.
\newblock {\em Toronto 2001, Calabi-Yau Varieties and Mirror Symmetry},
  hep-th/0402133:121--157.

\bibitem{dwork:def}
Bernard Dwork.
\newblock A deformation theory for the zeta function of a hypersurface.
\newblock {\em Proceedings of the International Congreee of Maths.}, pages
  247--259, 1962.

\bibitem{dwork:padic}
Bernard Dwork.
\newblock $p$-adic cycles.
\newblock {\em Pub. math. de l'I.H.\'{E}.S.}, 37:27--115, 1969.

\bibitem{kedlaya}
Kiran~Kedlaya et. al.
\newblock {\em $p$-adic Geometry: Lectures from the 2007 Arizona Winter School,
  University Lecture Series}, volume~45.
\newblock A.M.S., 2008.

\bibitem{gfunctions}
Bernard Dwork;~Giovanni Gerotto; and Francis~J. Sullivan.
\newblock {\em An Introduction to $G$-Functions}.
\newblock Princeton University Press, 1994.

\bibitem{griffiths}
Philip~A. Griffiths.
\newblock On the periods of certain rational integrals: {I}.
\newblock {\em The Annals of Mathematics}, 90--3:460--495, 1969.


\bibitem{ince}
E.L. Ince.
\newblock {\em Ordinary Differential Equations}.
\newblock Dover Publications, 1944.

\bibitem{katz:period}
Nicholas~M. Katz.
\newblock On the differential equations satisfied by period matrices.
\newblock {\em Publ. Math. I.H.E.S.}, 35:71--106, 1968.

\bibitem{katz2}
Nicholas~M. Katz.
\newblock On the intersection matrix of a hypersurface.
\newblock {\em Ann. Sci. \'{E}cole Norm. Sup.}, 4--2:583--598, 1969.

\bibitem{katz:esde}
Nicholas~M. Katz.
\newblock {\em Exponential Sums and Differential Equations}.
\newblock Princeton University Press, 1990.

\bibitem{katz:dwork}
Nicholas~M. Katz.
\newblock Another look at the {D}work family.
\newblock {\em Manin Festschrift}, to appear.


\bibitem{katz:oda}
N.M. Katz;~T. Oda.
\newblock On the differentiation of de {R}ham cohomology classes with respect
  to parameters.
\newblock {\em J. Math. Kyoto Univ.}, 8:199--213, 1968.

\bibitem{kloost}
	Remke~Kloosterman.
	\newblock The zeta function of monomial deformations of {F}ermat hypersurfaces.
	\newblock{\em Algebra \& Number Theory}, 1:421--450, 2007.
	
\bibitem{PARI2}
    PARI/GP, version {\tt 2.3.3}, Bordeaux, 2011, {http://pari.math.u-bordeaux.fr/}.


\bibitem{frv:ent}
Fernando Rodr\'{i}guez-Villegas.
\newblock {\em Experimental Number Theory}.
\newblock Oxford University Press, 2007.

\bibitem{rojaswan}
A.~Rojas-Leon and D.~Wan.
\newblock Moment zeta functions for toric {C}alabi-{Y}au hypersufaces.
\newblock {\em Comm. in Number Theory and Physics}, 1-3:539--578, 2007.

\bibitem{salerno}
	Adriana~Salerno.
	\newblock {\em Hypergeometric Functions in Arithmetic Geometry}.
	\newblock Thesis, University of Texas at Austin, 2009.
	

\bibitem{slater}
Lucy~Joan Slater.
\newblock {\em Generalized Hypergeometric Functions}.
\newblock Cambridge University Press, 1966.



	




\end{thebibliography}

\end{document}